\documentclass[11pt]{article}
\pdfoutput=1 
\usepackage{amssymb}
\usepackage[margin=1in]{geometry}

\usepackage{lineno}
\usepackage{graphicx}

\usepackage{moreverb,url}
\usepackage{makecell}
\usepackage[breaklinks=true, colorlinks=true, linkcolor={blue!90!black}, citecolor={blue!90!black}, urlcolor={blue!90!black}]{hyperref}

\newcommand\BibTeX{{\rmfamily B\kern-.05em \textsc{i\kern-.025em b}\kern-.08em
T\kern-.1667em\lower.7ex\hbox{E}\kern-.125emX}}

\usepackage[algo2e,linesnumbered,vlined,algoruled]{algorithm2e}

\usepackage{placeins}
\SetCommentSty{mycommfont}

\usepackage{geometry,xcolor,tabularx,nccmath}
\usepackage{amsmath,amsthm,amssymb,enumitem}
\usepackage{booktabs,caption,subcaption,mathtools,multirow}
\usepackage{graphicx,tabulary}
\usepackage{rotating}

\usepackage[utf8]{inputenc}
\usepackage{amsfonts}
\usepackage{fullpage}

\newcommand{\f}{f_{u^*,v^*,i^*,j^*}}

\usepackage[colorinlistoftodos]{todonotes}

\usepackage{placeins}
\usepackage{etoolbox}

\usepackage[nameinlink]{cleveref}
\hypersetup{colorlinks={true},allcolors={blue}}

\crefformat{section}{\S#2#1#3}
\crefformat{subsection}{\S#2#1#3}
\crefformat{equation}{#2Equation~\textup{(#1)}#3}
\crefformat{cons}{#2\textbf{Constraint}~\textup{\textbf{(#1)}}#3}
\labelcrefformat{cons}{\textup{#2(#1)#3}}
\crefformat{consset}{#2\textbf{Constraints}~\textup{\textbf{(#1)}}#3} 
\labelcrefformat{consset}{\textup{\textbf{#2(#1)#3}}}
\crefformat{func}{#2Function~\textup{(#1)}#3}
\labelcrefformat{func}{#2Function~\textup{(#1)}#3}
\crefformat{algo}{#2\textbf{Algorithm}~\textup{\textbf{#1}}#3}
\crefformat{objfunc}{#2\textbf{Objective function}~\textup{\textbf{(#1)}}#3}
\labelcrefformat{objfunc}{\textup{(#2#1#3)}}
\crefformat{ch}{#2Chapter~\textup{#1}#3}
\crefformat{sec}{#2Section~\textup{#1}#3}
\crefformat{tab}{#2Table~\textup{#1}#3}
\crefformat{fig}{#2Figure~\textup{#1}#3}
\crefformat{eq}{#2equation~\textup{(#1)}#3}
\labelcrefformat{eq}{\textup{(#2#1#3)}}
\crefformat{theo}{#2Theorem~\textup{#1}#3}
\crefformat{lem}{#2Lemma~\textup{#1}#3}
\crefformat{app}{#2Appendix~\textup{#1}#3}
\crefformat{scenario}{#2Scenario~\textup{#1}#3}
 \AtBeginDocument{%
    \crefname{figure}{Figure}{figures}%
}

\expandafter\def\expandafter\UrlBreaks\expandafter{\UrlBreaks
  \do\a\do\b\do\c\do\d\do\e\do\f\do\g\do\h\do\i\do\j%
  \do\k\do\l\do\m\do\n\do\o\do\p\do\q\do\r\do\s\do\t%
  \do\u\do\v\do\w\do\x\do\y\do\z\do\A\do\B\do\C\do\D%
  \do\E\do\F\do\G\do\H\do\I\do\J\do\K\do\L\do\M\do\N%
  \do\O\do\P\do\Q\do\R\do\S\do\T\do\U\do\V\do\W\do\X%
  \do\Y\do\Z}

\usepackage{amsmath,cases,eqparbox,xparse}

\makeatletter
\NewDocumentCommand{\eqmathbox}{o O{c} m}{%
  \IfValueTF{#1}
    {\def\eqmathbox@##1##2{\eqmakebox[#1][#2]{$##1##2$}}}
    {\def\eqmathbox@##1##2{\eqmakebox{$##1##2$}}}
  \mathpalette\eqmathbox@{#3}
}
\makeatother

\makeatletter
\let\@msm@th@eqref\eqref
\renewcommand{\eqref}[1]{%
  \begingroup
  \leavevmode
  \color{blue}%
  \hypersetup{linkbordercolor=[named]{blue}}%
  \@msm@th@eqref{#1}%
  \endgroup
}
\makeatother

\usepackage[normalem]{ulem} 
 
\crefname{defi}{definition}{definitions}
\Crefname{defi}{Definition}{Definitions}
\crefname{lemma}{lemma}{lemmas}
\Crefname{lemma}{Lemma}{Lemmas}
\crefname{assumption}{assumption}{assumptions}
\Crefname{assumption}{Assumption}{Assumptions}

\usepackage{multicol}
\usepackage[super]{nth}
\usepackage{authblk}
\providecommand{\keywords}[1]{\noindent\textbf{Keywords:} #1}

\usepackage{natbib}

\begin{document}

\title{Impact of Transit on Mobility, Equity, and Economy in the Chicago Metropolitan Region}


\author[1]{Omer Verbas}
\author[1]{Taner Cokyasar}
\author[2]{Seamus Joyce-Johnson}
\author[3]{Scott Wainwright}
\author[3]{Maeve Coates}
\author[1]{Aymeric Rousseau}
\author[2]{Jim Aloisi}
\author[2]{Anson Stewart}
\author[1]{Joshua Auld}

\affil[1]{Argonne National Laboratory 9700 S Cass Ave, Lemont, IL 60439, USA}

\affil[2]{Massachusetts Institute of Technology, 77 Massachusetts Ave, Cambridge, MA 02139, USA}

\affil[3]{Chicago Transit Authority, 567 W Lake St, Chicago, IL 60661, USA}

\maketitle


\begin{abstract}
Transit is essential for urban transportation and achieving net-zero targets. In urban areas like the Chicago Metropolitan Region, transit enhances mobility and connects people, fostering a dynamic economy. To quantify the mobility and selected economic impacts of transit, we use a novel agent-based simulation model POLARIS to compare baseline service against a scenario in which transit is completely removed. The transit-removal scenario assumes higher car ownership and results in higher traffic congestion, numerous activity cancellations, and economic decline. In this scenario, average travel times increase by 14.2\% regionally and 34.7\% within the City of Chicago. The resulting congestion causes significant activity cancellations despite increased car ownership: 11.8\% of non-work and 2.8\% of work/school activities regionally, totaling an 8.6\% overall cancellation rate. In the city, non-work cancellations would reach 26.9\%, and work/school cancellations 7.3\%, leading to a 19.9\% overall cancellation rate. The impact varies between groups. Women and lower-income individuals are more likely to cancel activities than men and higher-income groups. Women account for 53.7\% of non-work and 53.0\% of total cancellations. The lowest 40\% income group experiences 50.2\% of non-work and 48.0\% of overall cancellations. Combined, activity cancellations, travel time losses, and increased car ownership cost the region \$35.4 billion. With annual public transit funding at \$2.7 billion, the ratio is 13 to 1, underscoring transit's critical role in mobility, equity, and economic health. 
\end{abstract}


\section{Introduction}
The global population is gradually rising and is projected to reach 10.2 billion by 2060 \citep{morse_global_nodate}. In 2022, car ownership in the U.S. reached a peak, with 91.7\% of households having at least one vehicle \citep{valentine_car_2023}. Nevertheless, public transit plays a crucial role in addressing various transportation issues, including reducing congestion, minimizing environmental impact, and mitigating economic effects.  

Since the COVID-19 pandemic, U.S. transit agencies have faced a historic decline in ridership and revenue. Although nationwide ridership has rebounded to over 70 percent of pre-pandemic levels \citep{dickens_public_2023}, the downturn prompted discussions about potentially eliminating underperforming bus and rail routes or reducing service frequencies. While removing poorly performing routes can enhance overall traffic flow when paired with improved service on the remaining routes \citep{ewing_travel_2010}, it may also result in increased congestion if these routes are not replaced with viable alternatives, or if the remaining transit options are inefficient \citep{buehler_demand_2012}. It is important to note that reductions in service can have severe environmental consequences \citep{griswold_unintended_2014}. 

In response to global net-zero environmental targets, regulators are developing innovative strategies to meet these goals. Electrifying the whole transportation sector is a highly favored target and will be crucial in achieving net-zero emissions, but enhancing its efficiency by promoting shared transportation modes, such as mass transit, on-demand transit, car-pooling, and integrating transit with active modes (e.g., shared bikes and scooters), is also important. Encouraging the use of shared transportation modes can significantly contribute to reducing overall emissions. Studies showed that well-implemented transportation policies can effectively reduce greenhouse gas emissions \citep{cheng_minimizing_2018,chester_time-based_2016, griswold_tradeoffs_2013}.  

In the U.S., low-density land use patterns and widespread reliance on private vehicles often result in an increase in car ownership and, consequently, vehicle miles traveled (VMT) and traffic congestion when transit services are removed. However, the full extent of these impacts and their implications for urban mobility, equity, and economic factors has not been extensively studied in the literature. To address this gap, this study analyzes a scenario for the Chicago Metropolitan Region in which all transit services are removed. We use POLARIS (Planning and Operations Language for Agent-based Regional Integrated Simulation) to model \citep{auld_polaris:_2016} and analyze the implications of such a removal on various dimensions, including mobility, equity, and economic performance. POLARIS allows for a detailed simulation of how transit service reductions would alter travel behavior, affect different socio-economic groups, and change regional economic conditions. Through this analysis, we aim to provide a comprehensive understanding of the potential consequences and inform policy decisions to better manage and plan for changes in transit services. 

The main contribution of this study is to examine the value of transit services in the Chicago Metropolitan Region, in terms of mobility, equity, and the economy. By considering a total removal of transit, the study predicts significant consequences of an extreme scenario. While it may not be realistic to imagine a complete removal of transit all at once, the reduction in the level of service and the resulting loss of boardings, hence, revenue could cause a downward spiral of transit eventually being completely removed. Furthermore, modeling a hypothetical travel regime without transit allows us to understand the full force of the societal benefits afforded by transit service, which do not come across in models of incremental service changes. Utilizing the agent-based POLARIS framework, this research captures the downstream effects of transit removal on modal shifts, congestion, activity delays and cancellations, and their impact on equity and the economy. Without a comprehensive and detailed agent-based modeling framework like this one, it is challenging to capture the complex interactions between a major change in the service level of one mode and the resulting impacts on traffic and the decision-making processes of millions of agents in a metropolitan region. Furthermore, a detailed breakdown by agent allows for a closer investigation of how impacts are distributed across different demographic groups. 

\section{Literature review} \label[sec]{lit_rev}

In the existing academic literature, there are limited studies addressing the complete removal of transit services, with only two such studies providing a qualitative analysis of the issue \citep{nguyen-phuoc_how_2018,nguyen-phuoc_transit_2018}. These studies approached the problem from a qualitative perspective rather than a simulation-based analytical approach adopted in this study. Such studies are found in industry; for example, agencies and non-governmental organizations in Boston and the Washington, D.C. region have used a transit removal modeling approach to convey the full value of public transit service in those regions \citep{aecom_cost_2013}. Much of the literature has concentrated on temporary disruptions in transit services \citep{adele_exploring_2019,fearnley_competition_2018,nguyen-phuoc_exploring_2018}. Additionally, there is a body of work exploring various strategies for reducing transit budget expenditures, including the removal of low-demand bus stops and routes \citep{chien_genetic_2001,mei_planning_2021,torabi_limited-stop_2019,wagner_benefitcost_2014}, reductions in service frequency \citep{verbas_optimal_2013,verbas_exploring_2015}, network redesign \citep{bertsimas_data-driven_2021,ng_redesigning_2024}, and privatization efforts \citep{giuliano_estimating_1986,teal_privatization_1986}. These studies offer valuable insights into different aspects of transit service management but do not fully address the impacts of complete service removal in a simulation context. 

Travel behavior of users under a complete transit service cut was examined in \citet{nguyen-phuoc_how_2018}. The qualitative study highlighted various adaptive strategies, including increased reliance on private vehicles, carpooling, and changes in travel patterns. The research highlighted that users' responses are influenced by factors such as trip purpose, availability of alternatives, and personal circumstances. Similar to the findings of this study, the authors stated that riders either switch to alternative modes such as cars, bikes, walking or cancel their trips \citep{nguyen-phuoc_transit_2018}. 

In \citet{adele_exploring_2019}, the authors focused on the behavior of suburban train users during service disruptions. The study revealed that users adapt by switching to alternative transport modes or altering departure times to avoid peak disruptions. The competition and substitution between different public transport modes was analyzed in \citet{fearnley_competition_2018}. While some users opt for alternative public transport modes, others may turn to private vehicles or ride-sharing services. 

Public transit plays a key role in enabling daily activities beyond the traditional commute. Commutes to and from work account for less than 20\% of trips nationwide \citep{mcguckin_summary_2018}. The remaining trips constitute everything from shopping to running errands to social activities. A subset of these is known as \emph{mobility of care} trips, and they are typically essential for the overall quality of life of most people. Mobility of care in the transportation context focuses on how people travel to activities classified as care work, a term defined as \emph{the unpaid labor performed by adults for children and other dependents, including labor related to the upkeep of a household} \citep{madariaga_mobility_2013}. 

A review of Chicago Metropolitan Agency for Planning (CMAP) survey data on trip purposes \citep{chicago_metropolitan_agency_for_planning_household_2019} revealed that women accounted for 61\% of care trips and that lower-income residents made a greater proportion of care trips (23\%) than upper income residents (20\%) in the Chicago region. Seven of the 28 trip purposes were flagged as relating to mobility of care: shopping, drive-thru errands, health care visit for someone else, visiting a person staying at the hospital, non-shopping errands, drop-off/pick-up, and accompanying someone else. Performing a holistic evaluation of a transportation system requires accounting for these trips. 

In \citet{torabi_limited-stop_2019}, the authors focused on limited-stop bus services as a strategy to reduce unused capacity in transit networks. Their research highlighted how implementing limited-stop services can address capacity issues and improve the efficiency of transit operations. By analyzing the trade-offs between service frequency and capacity utilization, the study shows that limited-stop services can effectively manage excess capacity and enhance overall network performance. 

A benefit-cost evaluation method for transit stop removals was presented in \citet{wagner_benefitcost_2014}. The study provided a framework for assessing the impacts of removing transit stops on service efficiency and passenger access. By evaluating the trade-offs between operational benefits and accessibility impacts, the method helps transit agencies make informed decisions about stop removals. The study underscores the importance of balancing service improvements with potential negative effects on riders, particularly in underserved areas. 

The optimal allocation of service frequencies across transit network routes and time periods was explored in \citet{verbas_optimal_2013}. Their research addressed the challenges of determining the ideal frequency of service to maximize network performance and meet varying demand levels. Using bus route patterns, they proposed a formulation and solution method for frequency allocation, highlighting the benefits of tailored frequency adjustments for different times and routes. In an extended work \citep{verbas_exploring_2015}, the authors investigated the trade-offs involved in frequency allocation using large-scale urban systems. Their study develops a methodology for exploring these trade-offs, emphasizing the need to balance service frequency with operational costs and demand coverage. The research demonstrated that optimizing frequency allocation can lead to significant improvements in transit service quality and efficiency. 

In a recent study, the redesign of large-scale multimodal transit networks with the integration of shared autonomous mobility services was investigated \citep{ng_redesigning_2024}. The study explored how autonomous vehicles can be incorporated into existing transit networks to enhance connectivity and efficiency. By leveraging shared autonomous mobility services, transit networks can improve service coverage and operational flexibility. The study provided insights into the potential benefits and challenges of integrating these technologies into large-scale transit systems, emphasizing the transformative impact of autonomous mobility on network design. 

\section{Methodology}\label[sec]{methodology}
This research utilizes POLARIS \citep{auld_polaris:_2016}, a comprehensive and large-scale software that integrates multi-agent, activity-based travel demand, and multimodal traffic and transit assignment. It models the movements of both individuals and freight over a 24-hour period. 

\subsection{Activity-Based Demand Modeling}
POLARIS synthesizes a population of households and individuals using data from U.S. Census tracts, Public Use Microdata Areas (PUMAs), and the American Community Survey (ACS). This synthesis includes choices for home, school, and work locations \citep{auld_efficient_2010}. Based on the socioeconomic details at the personal and household levels, activities are generated for each individual over a 24-hour period. The start time and duration of each activity are determined using a hazard-based model \citep{auld_dynamic_2011}. A multinomial logit model determines the activity location choice, and a nested logit model determines the mode choice \citep{auld_activity_2012}. In addition to choices regarding mode, location, start time, and duration, the choice of travel party is also modeled \citep{auld_planning-constrained_2011}. Finally, a conflict resolution model is implemented to resolve scheduling conflicts both within an individual's plans and among individuals within the household \citep{auld_modeling_2009}. 

\subsection{Freight Modeling}
Similar to the passenger travel, POLARIS also synthesizes firms and establishments using the CRISTAL model \citep{stinson_introducing_2022}. These firms and establishments ultimately generate truck and e-commerce delivery trips \citep{cokyasar2023optimization,cokyasar2023comparing,cokyasar_time-constrained_2023,davatgari2024electric,subramanyam2022joint,zuniga-garcia_freight_2023}. While the details of these models are beyond the scope of this study, it is important to emphasize that transportation is a complex-adaptive system where the decision of each agent affects the decision making of other groups. As an example, if there are more trucks on arterial roads, the buses using these roads are slowed, which in turn causes a mode shift to passenger cars, that generates more traffic, which would ultimately also affect the speed of the trucks. 

\subsection{Multimodal Supply Modeling}
Once the population is synthesized and the activity schedules of each individual including timing, location, mode, and party are generated; and all the truck and delivery trips are created utilizing the CRISTAL module, the simulation kicks off. Each trip is routed using a time-dependent intermodal A* algorithm \citep{verbas_time-dependent_2018}. This algorithm finds routes for driving, taxi/transportation network company (TNC) and truck/delivery modes on the roadway network, and walk and bike modes on the walking and biking networks, respectively. Intermodal routes, such as walk-to-transit, drive-to-transit, and TNC as a first-mile-last-mile (FMLM) to/from transit, are generated utilizing the required combination of active, transit, and roadway networks brought together in a multimodal network representation. The network equilibrium is established within a dynamic transit and traffic assignment framework using a framework of information-mixing between prevailing and historical network conditions \citep{auld_agent-based_2019}.  

The traffic is simulated using a multi-class traffic flow model in Lagrangian coordinates \citep{de_souza_polaris-lc_2024}, where every class of vehicle has a different speed-spacing relationship. This is an extension on the mesoscopic traffic flow model \citep{souza_mesoscopic_2019} previously used in POLARIS, where every link used to have a single (average) speed-spacing relationship based on the traffic mix on the link. The transit buses move on the congestable roadway network, whereas the trains run on a non-congestable network adhering to the General Transit Feed Specification (GTFS) travel times \citep{noauthor_general_nodate}. The interaction between the transit vehicles and passengers such as boarding, seating, standing, and alighting are also simulated. 

POLARIS also simulates the movement and decision making of TNC agents for services of today including single rides, pooled rides, corner-to-corner trips, integration with transit as a FMLM service, on-demand transit services, as well as more centralized and autonomous business models of the potential future \citep{gurumurthy_integrating_2020,gurumurthy_system_2021,zuniga-garcia_integrating_2022,sarma_-demand_2023,lu_revitalizing_2023}. 

\subsection{Integration and Feedback Loops}
As stated in the previous subsection, the network equilibrium is established within a dynamic transit and traffic assignment framework. However, network equilibrium is only an inner loop to establish consistency between routed travel times at the beginning of a trip and experienced travel times resulting from simulating the traffic on a congestable multimodal network with multiple vehicles of different classes. On top of that, if a trip takes much longer than expected, an agent may decide to change their route or even their mode of travel. Moreover, a late arrival at an activity may trigger shortening of the said activity, or maintaining the activity duration and postponing/cancelling the subsequent activity, and so on. In summary, the countless interactions between travelers, freight agents, transit and TNC service providers are captured in a continuous fashion during one iteration of the simulation, and the overall performance of each mode and agent affects their decision-making at the beginning of the subsequent iteration. 

\subsection{Model Validation }
POLARIS has been validated and used in multiple large-scale studies. Interested readers are encouraged to review the following publications: A study evaluating new transportation technologies such as connectivity, automation, and sharing \citep{rousseau_smart_2020}, validating POLARIS ride-hail simulation through back-casting in Chicago \citep{verbas_validating_2021}, optimization of transit frequencies in Chicago \citep{verbas2024modeling}, and a large-scale evaluation of mobility, technology and demand scenarios in the Chicago region using POLARIS. 

\section{Scenarios}
In this study, a baseline network of 2025 is considered for the Chicago Metropolitan Region. There are four major transit agencies serving in the region: 

\begin{itemize}
\item Chicago Transit Authority (CTA) 
    \begin{itemize}
        \item Urban bus – 123 routes
        \item Urban rail (metro) – 8 routes 
    \end{itemize}
\item Pace 
    \begin{itemize}
        \item Suburban bus – 134 routes 
    \end{itemize}
    
\item Metropolitan Rail (METRA) 
    \begin{itemize}
        \item Commuter rail – 11 routes 
    \end{itemize}

\item Northern Indiana Commuter Transportation District (NICTD) 
    \begin{itemize}
        \item South Shore Line (SSL), a commuter rail line between South Bend, IN and downtown Chicago – 1 route 
    \end{itemize}
\end{itemize}

\begin{figure}[!ht]
    \centering
    \includegraphics[width=1\linewidth]{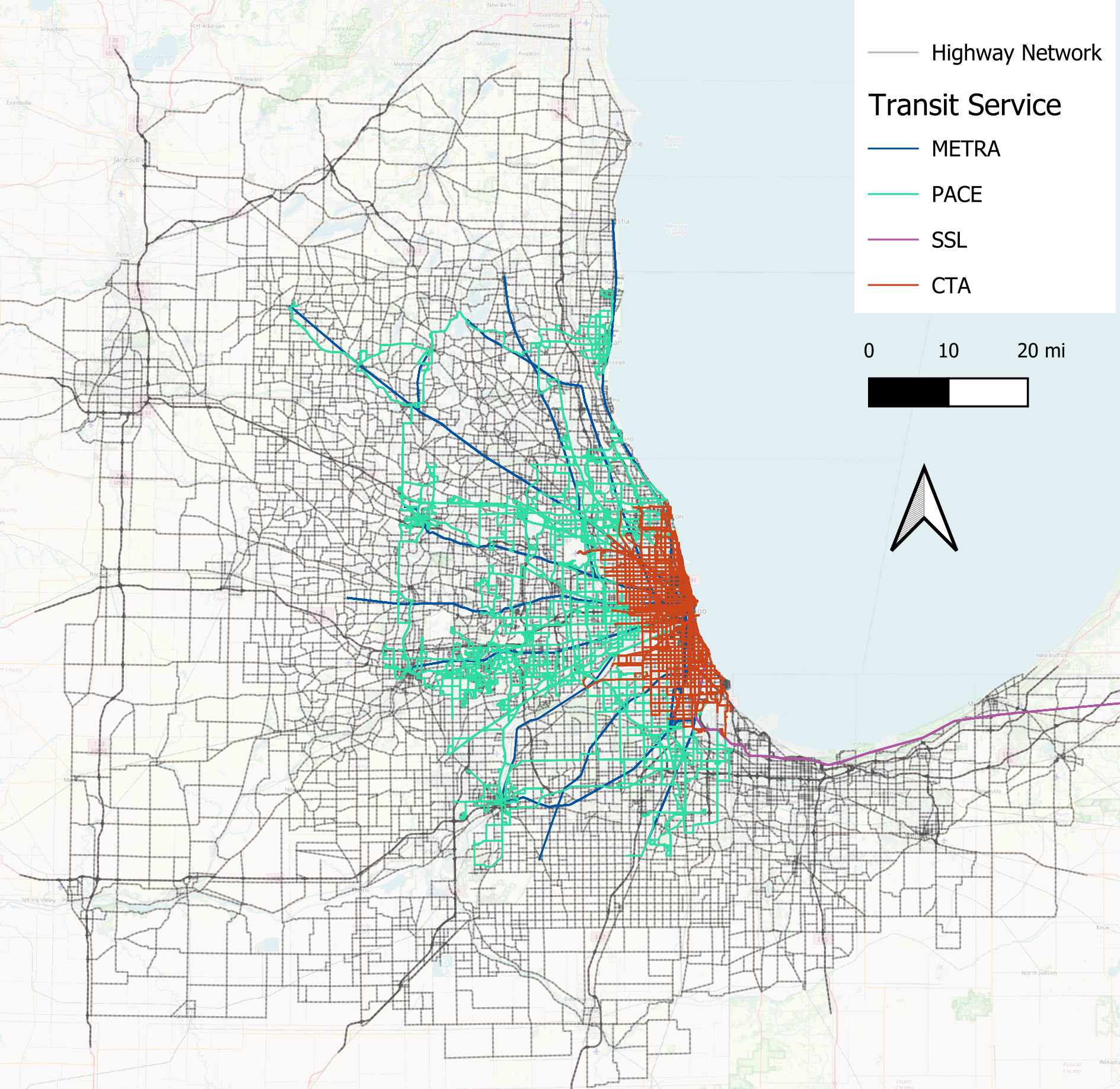}
    \caption{Chicago Metropolitan Region transportation network map.}
    \label[fig]{layout}
\end{figure}

On the roadway side there are 84,948 links and 35,583 nodes that are modeled. See \cref{layout} for the multimodal network representation of the Chicago Metropolitan Region. Two scenarios are defined: 
\begin{itemize}
    \item 2025 baseline 
    \item Transit removal with increased car ownership. 
\end{itemize}

\begin{figure}[!ht]
    \centering
    \includegraphics[width=1\linewidth]{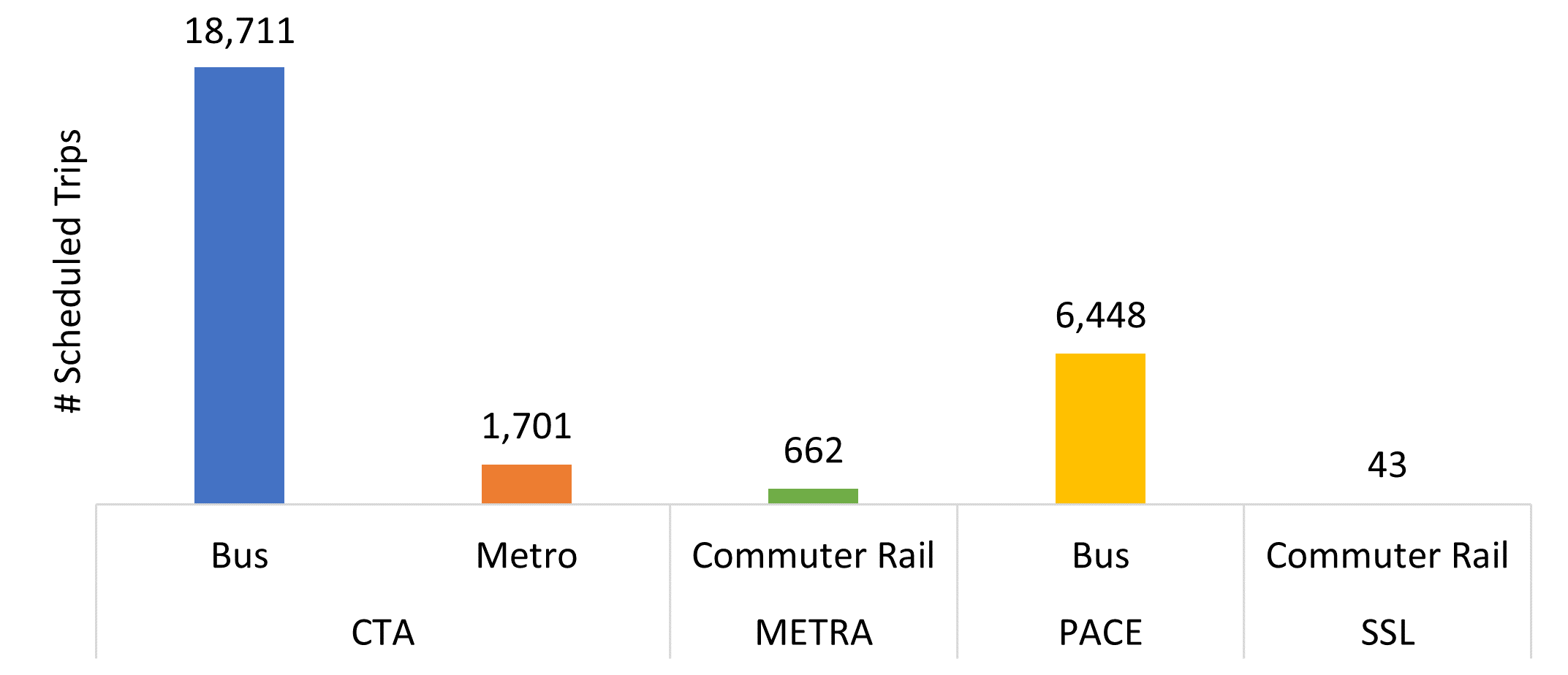}
    \caption{Number of scheduled transit revenue trips by agency and mode under the 2025 baseline scenario.}
    \label[fig]{trips}
\end{figure}

For the baseline scenario, the number of scheduled transit revenue trips is shown in \cref{trips}. Regional Transportation Authority (RTA), a state government agency that coordinates the Chicago region’s transit system \citep{noauthor_regional_nodate}, has shared that they do not expect a major service recovery to pre-COVID levels for Pace and METRA. As a result, Tuesday, March \nth{19}, 2024 is selected from their General Transit Feed Specification databases \citep{noauthor_general_nodate,noauthor_route_nodate,noauthor_developers_nodate}. On the other hand, CTA expects to go back to its full pre-COVID level of service by 2025. As a result, for the CTA and SSL, Tuesday, October \nth{8}, 2019 is selected from the GTFS \citep{noauthor_gtfsscheduled_nodate,noauthor_south_nodate}.  

\begin{figure}[!ht]
    \centering
    \includegraphics[width=1\linewidth]{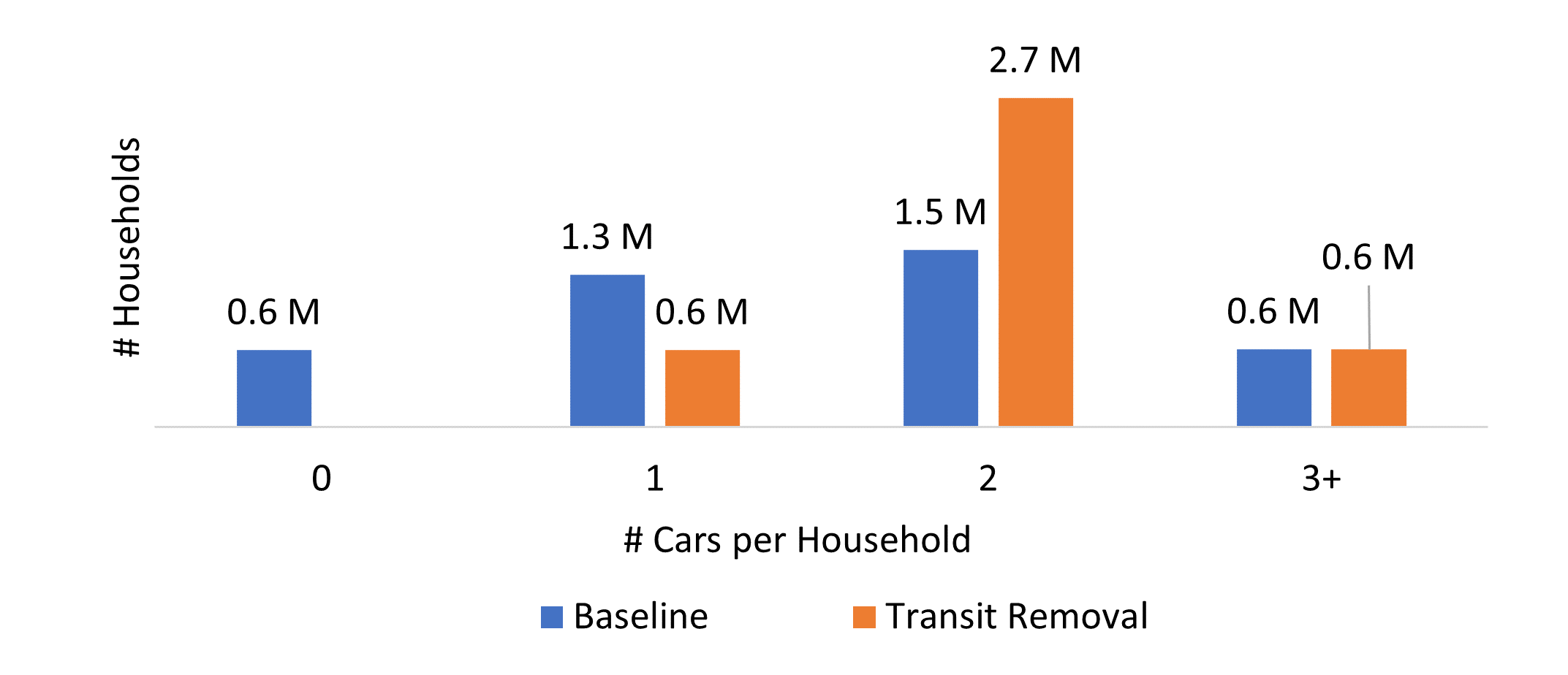}
    \caption{Household vehicle ownership distribution under both scenarios.}
    \label[fig]{vehicle_ownership}
\end{figure}

For the transit removal scenario, all agencies, their routes, and their scheduled trips are removed from the multimodal network. Moreover, it is assumed that all households with no cars buy a car. Households with one car buy an additional car. These assumptions result in a total increase of 1.9 million (30\%) cars owned in the region. See \cref{vehicle_ownership} for the car ownership distribution under both scenarios. 

\section{Results}
This section summarizes the results under the categories of boardings and mode share, congestion, activity cancellations, and economic impact. 

\subsection{Boardings and Mode Share}

\begin{figure}[!ht]
    \centering
    \includegraphics[width=1\linewidth]{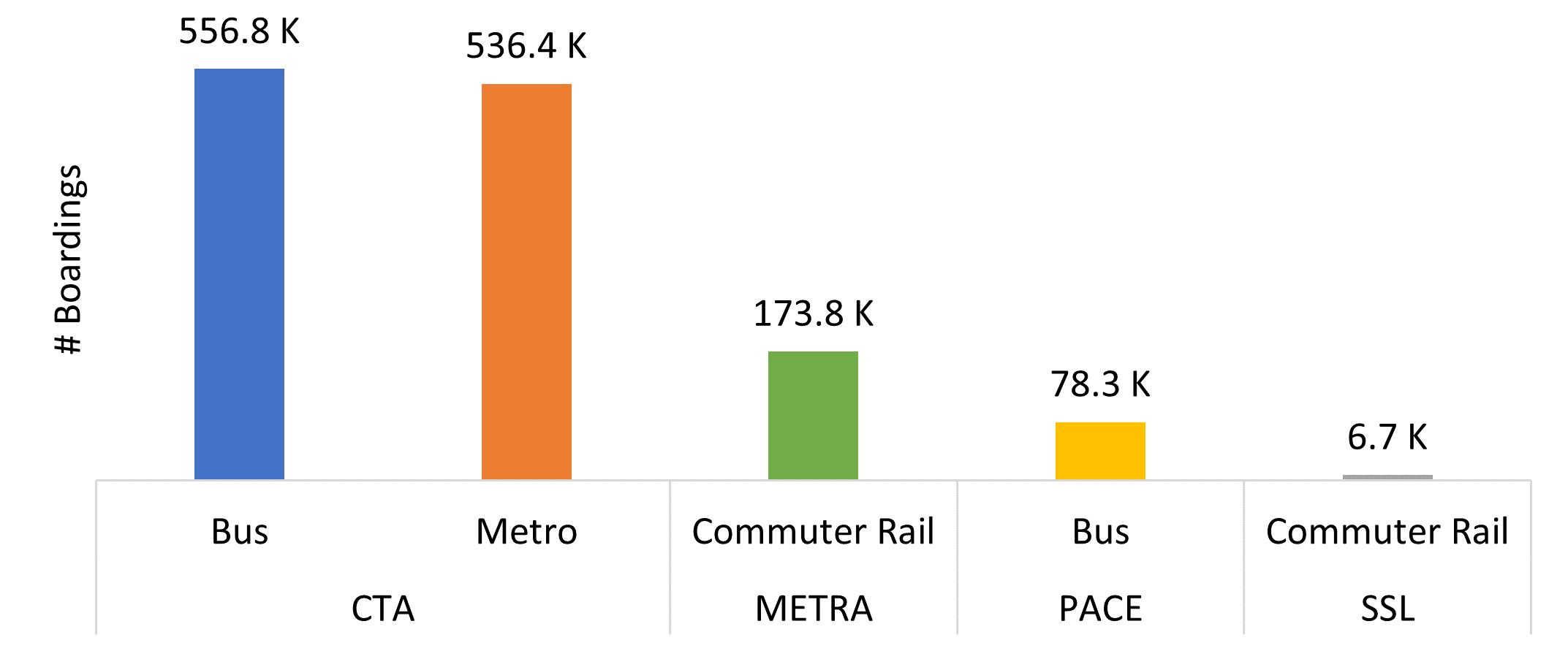}
    \caption{Transit boardings by agency and mode under the 2025 scenario.}
    \label[fig]{boardings}
\end{figure}

In the 2025 baseline, there are 1.35M transit boardings in the metropolitan region, which is substantially (30.2\%) below the 1.93M in October 2019 \citep{verbas_validating_2021}. Boardings drop to zero in the transit removal scenario. See \cref{boardings} for the transit boardings under the 2025 baseline scenario. 

\begin{figure}[!ht]
    \centering
    \includegraphics[width=1\linewidth]{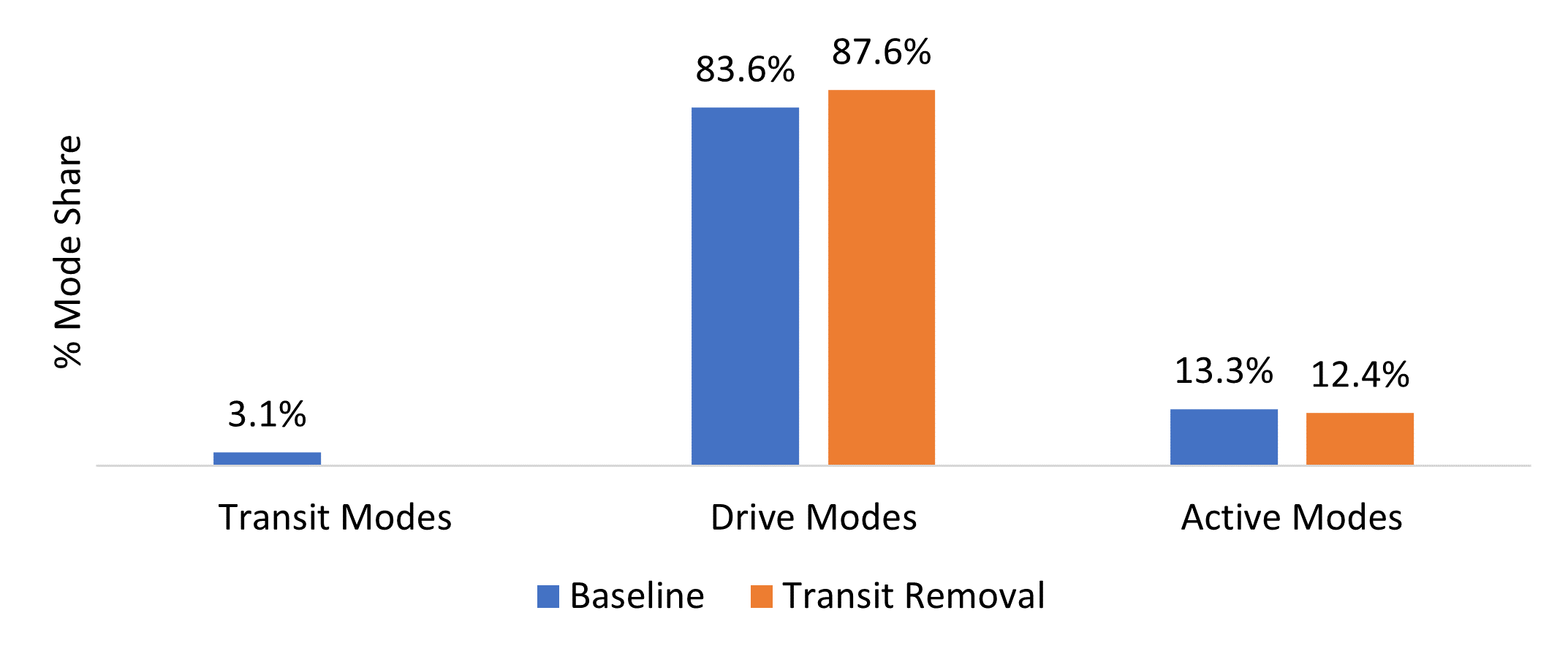}
    \caption{Mode share under both scenarios.}
    \label[fig]{mode_share}
\end{figure}

\cref{mode_share} shows that the mode share of transit modes – walk to transit, drive to transit, FMLM to/from transit – is 3.1\% under the 2025 scenario, which is already a massive (42.2\%) drop from the 5.4\% from the October 2019 baseline \citep{verbas_validating_2021}. While 1.35M boardings or 3.1\% transit mode share do not seem to account for so many trips, it has an outsize impact on traffic, which will be demonstrated in the upcoming section. 

\subsection{Congestion}
\begin{figure}[!ht]
    \centering
    \includegraphics[width=1\linewidth]{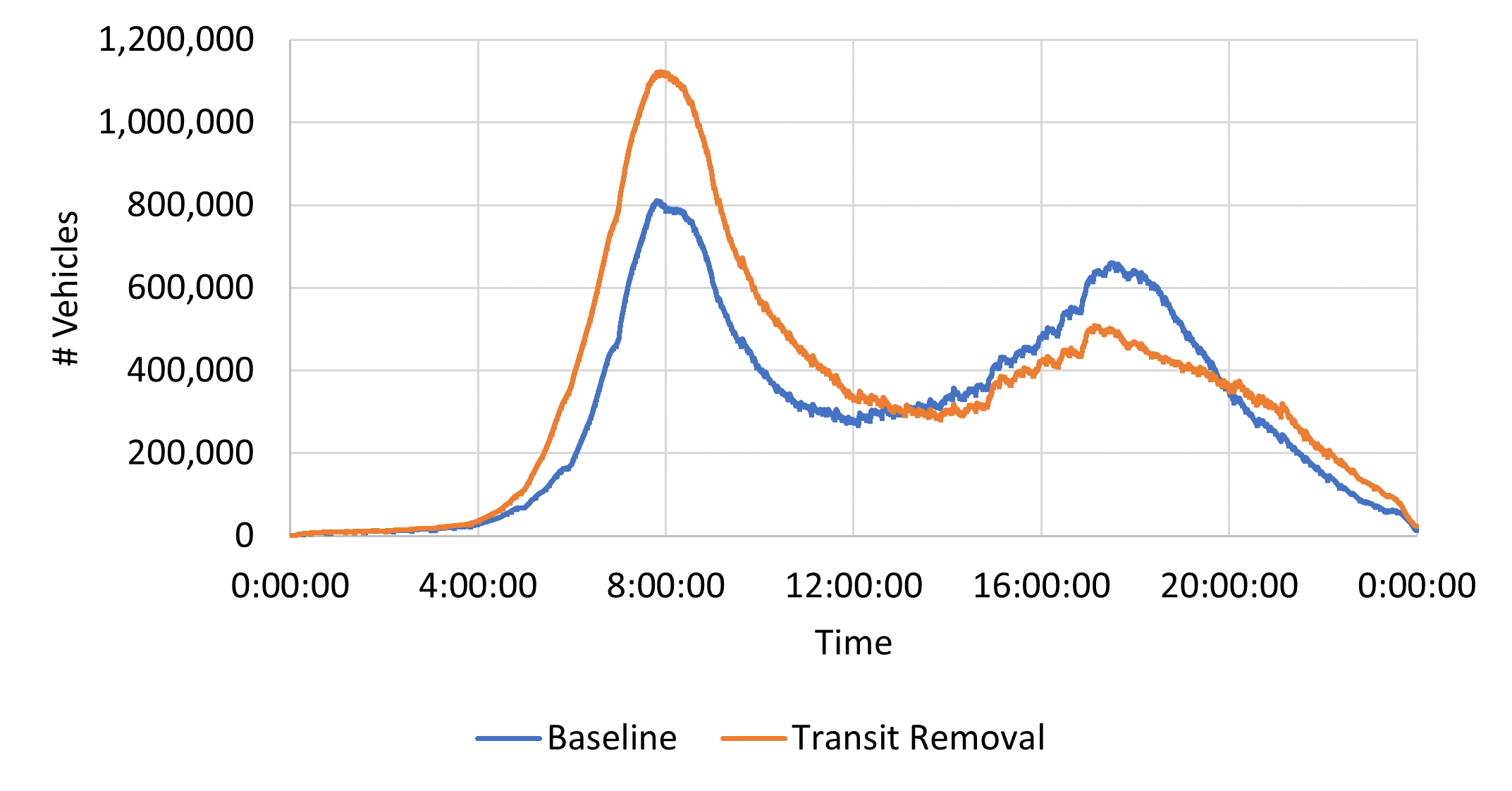}
    \caption{Number of roadway vehicles in the network at a given time under both scenarios.}
    \label[fig]{vehicles}
\end{figure}

The removal of transit services causes both activity cancellations, which will be covered in the following section, and substantial congestion. See \cref{vehicles} for the number of roadway vehicles in the network. Note that due to massive congestion in the AM peak under the transit removal scenario, many activities get cancelled, which results in some congestion relief for the rest of the day.  

\begin{figure}[!ht]
    \centering
    \includegraphics[width=1\linewidth]{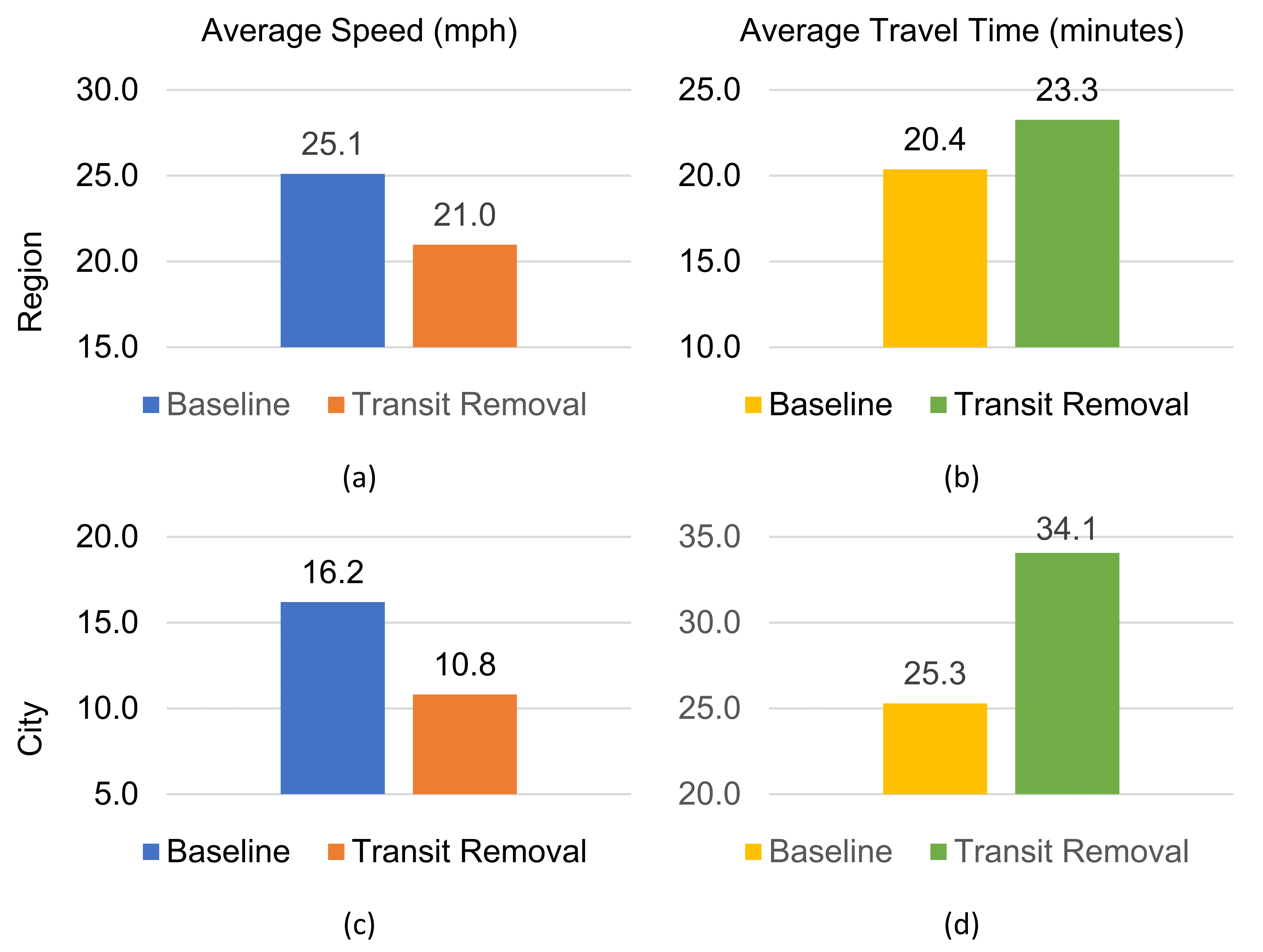}
    \caption{Average speeds and travel times in the region and in the city under both scenarios.}
    \label[fig]{average_speeds}
\end{figure}

\cref{average_speeds} demonstrates the change in average speeds and travel times for the region and for the city under both scenarios. Due to removal of transit and increase in car ownership, speeds decrease by 16.5\% (\cref{average_speeds}a), and travel times increase by 14.2\% (\cref{average_speeds}b) across the region. The impact of transit removal is exacerbated in the city. Speeds decrease by 33.2\% (\cref{average_speeds}c), and travel times increase by 34.7\% (\cref{average_speeds}d). 

\subsection{Activity Cancellations}
Due to the combined effect of increased congestion and transit service removal, a substantial portion of activities are cancelled. See \cref{activities_region} for the number of activities by activity type and scenario, as well as percentage changes in the entire metropolitan region. Around 1.7 million (11.8\%) non-work activities are cancelled. The work/school activity cancellation rate is at around 200 thousand (2.8\%) per day. The overall daily activity cancellation rate is at around 1.9 million (8.6\%).

\begin{table}[!ht]
  \footnotesize
  \caption{Number of activities by type and scenario in the region.}\label[tab]{activities_region}
  \begin{center}
    \begin{tabular*}\textwidth{l@{\extracolsep{\fill}}lll}
      \toprule
Activity type   & Baseline & Transit removal & Percentage change \\
\midrule
Eat out         & 2.15 M   & 1.93 M          & -10.1             \\
Errands         & 1.28 M   & 1.22 M          & -4.8              \\
EV charging     & 0.3 K    & 0.4 K           & 28.8              \\
Healthcare      & 893.7 K  & 845.2 K         & -5.4              \\
Leisure         & 1.67 M   & 1.42 M          & -14.8             \\
Part-time work  & 650.9 K  & 629.4 K         & -3.3              \\
Personal        & 327.4 K  & 304.4 K         & -7.0              \\
Pickup-dropoff  & 2.30 M   & 1.73 M          & -24.7             \\
Religious-civic & 386.1 K  & 330.7 K         & -14.4             \\
School          & 1.91 M   & 1.90 M          & -0.5              \\
Service         & 509.3 K  & 473.5 K         & -7.0              \\
Shop-major      & 666.5 K  & 624.0 K         & -6.4              \\
Shop-other      & 2.79 M   & 2.58 M          & -7.3              \\
Social          & 1.18 M   & 1.01 M          & -14.2             \\
Work            & 4.23 M   & 4.06 M          & -4.0              \\
Work at home    & 862.9 K  & 848.6 K         & -1.7              \\
\midrule
Total           & 21.80 M  & 19.92 M         & -8.6              \\
      \bottomrule
    \end{tabular*}
  \end{center}
\end{table}

These activity cancellations span across categories: with the exception of the most essential activities like school and work trips, large drops were observed in key categories like healthcare and pick-up and drop-off. Furthermore, activities that translate to consumer spending—such as eating out, errands, shopping, and leisure—showed the highest rates of activity cancellation. 

Activity cancellations are more severe in the city as result of even higher congestion increase and higher transit service loss. See \cref{activities_city} for the number of activities by activity type and scenario, as well as percentage changes in the City of Chicago. Around 1.0 million (26.9\%) non-work activities are cancelled. The work/school activity cancellation is at around 150 thousand (7.3\%). The overall activity cancellation is at around 1.2 million (19.9\%). 

\begin{table}[!ht]
  \footnotesize
  \caption{Number of activities by type and scenario in the city.}\label[tab]{activities_city}
  \begin{center}
    \begin{tabular*}\textwidth{l@{\extracolsep{\fill}}lll}
      \toprule
Activity type   & Baseline & Transit removal & Percentage change \\
\midrule
Eat out         & 581.3 K  & 422.3 K         & -27.4             \\
Errands         & 336.9 K  & 272.9 K         & -19.0             \\
EV charging     & 0.1 K    & 0.1 K           & 47.4              \\
Healthcare      & 231.6 K  & 187.4 K         & -19.1             \\
Leisure         & 438.6 K  & 296.0 K         & -32.5             \\
Part-time work  & 155.5 K  & 141.3 K         & -9.1              \\
Personal        & 86.4 K   & 67.3 K          & -22.1             \\
Pickup-dropoff  & 578.2 K  & 382.6 K         & -33.8             \\
Religious-civic & 104.1 K  & 72.3 K          & -30.5             \\
School          & 490.9 K  & 485.1 K         & -1.2              \\
Service         & 132.4 K  & 103.2 K         & -22.1             \\
Shop-major      & 176.6 K  & 136.4 K         & -22.7             \\
Shop-other      & 750.0 K  & 574.8 K         & -23.4             \\
Social          & 313.7 K  & 211.1 K         & -32.7             \\
Work            & 1.15 M   & 1.03 M          & -10.3             \\
Work at home    & 253.2 K  & 242.6 K         & -4.2              \\
\midrule
Total           & 5.78 M   & 4.63 M          & -19.9             \\
      \bottomrule
    \end{tabular*}
  \end{center}
\end{table}

While the congestion impact of transit removal is severe, some of it is mitigated by cancellation of activities. However, transportation is a derived demand, where the main objective is to provide access to work, school, and essential and non-essential services. Taking away this access through congestion and removal of transit service has an outsize impact not only on mobility but also on equity and the economy. 

\subsection{Equity}

\begin{figure}[!ht]
    \centering
    \includegraphics[width=1\linewidth]{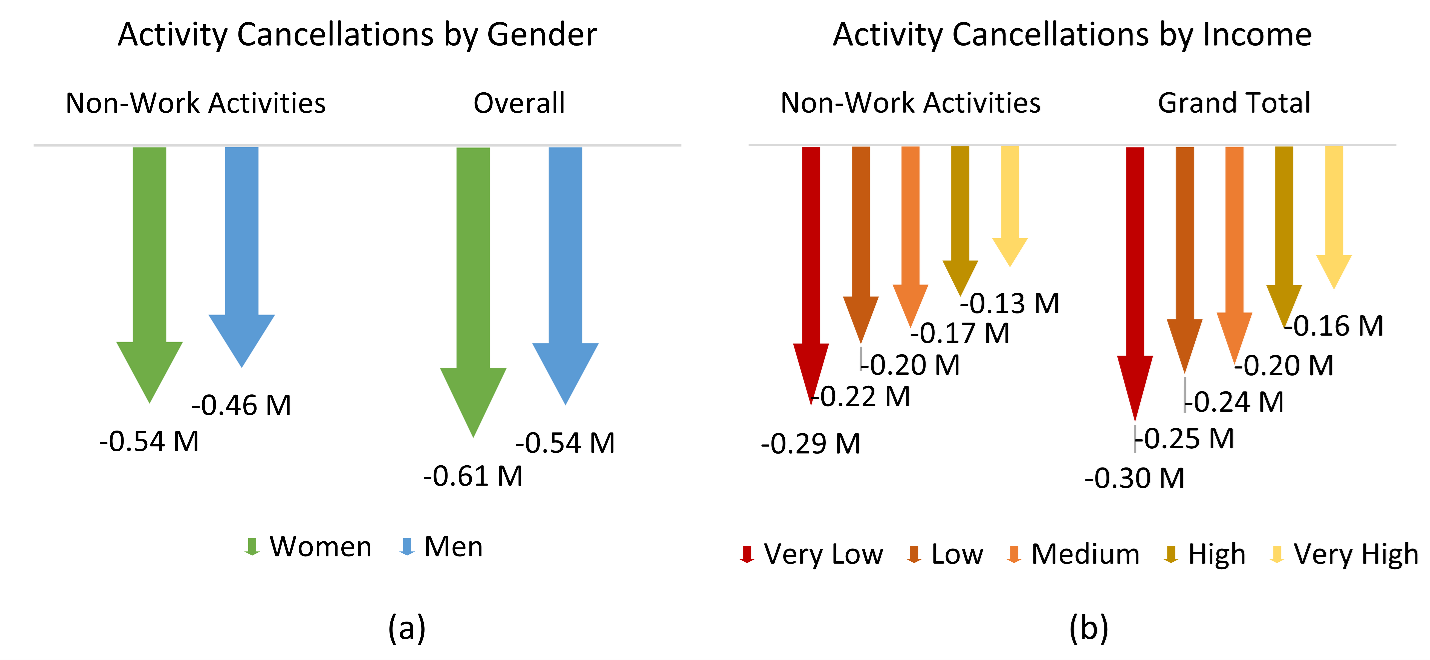}
    \caption{Activity cancellations in the City of Chicago by a) gender and b) income.}
    \label[fig]{activity_cancel}
\end{figure}

\cref{activity_cancel} demonstrates the number of activity cancellations by gender and income in the City of Chicago. Women, 51\% of the population \citep{census_reporter_census_2022}, account for 53.7\% of non-work activity cancellations and 53.0\% of overall activity cancellations (\cref{activity_cancel}a).

The households are categorized by income into five groups of equal size (quintiles). The lowest 20\% income group (very low) accounts for 28.5\% of non-work activity cancellations and 26.3\% of overall activity cancellations. The lowest 20\% and the next lowest 20\% (very low and low) income groups together account for 50.2\% of non-work activity cancellations, and 48.0\% of overall activity cancellations (\cref{activity_cancel}b). 

\begin{figure}[!ht]
    \centering
    \includegraphics[width=1\linewidth]{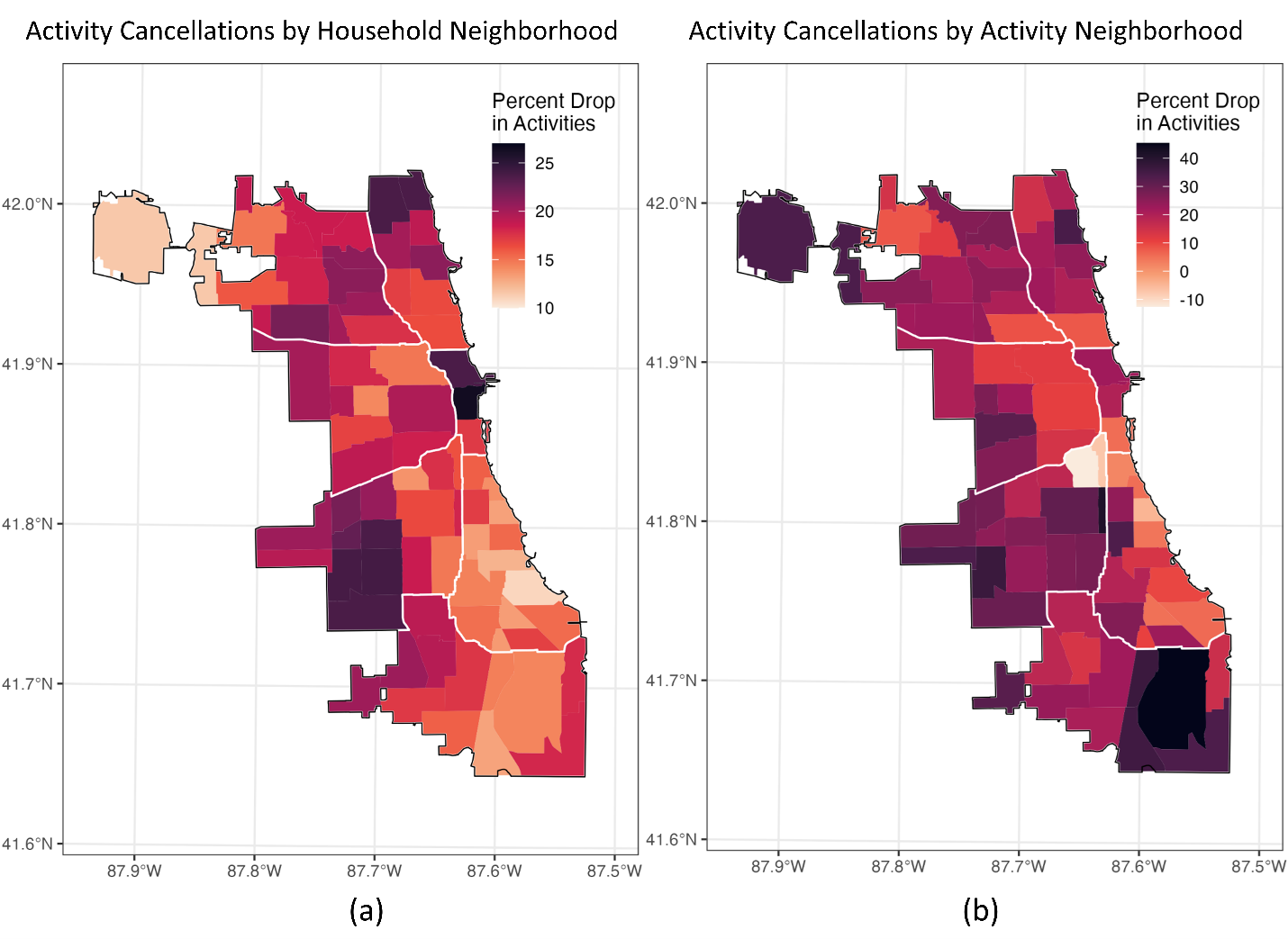}
    \caption{Activity cancellations in the City of Chicago by a) household neighborhood and b) activity neighborhood.}
    \label[fig]{activity_cancel2}
\end{figure}

Spatial visualization of activity cancellation across city neighborhoods reveals the uneven distribution of the benefits of transit across the Chicago (\cref{activity_cancel2}). Each activity modeled by POLARIS is associated with both the location the activity takes place, as well as the household location of the resident undertaking the activity. Households that cancel activities affect both their own quality of life and the neighborhoods where they would have gone (\cref{activity_cancel2}a). Cancellations by activity location concentrate in the outer neighborhoods of the city, threatening businesses and institutions (\cref{activity_cancel2}b). Seven of the ten neighborhoods with the highest activity decreases are on the Far South and Southwest Sides, indicating severe consequences for the economic vibrance of many of Chicago’s most disadvantaged neighborhoods. 

According to CMAP survey data \citep{chicago_metropolitan_agency_for_planning_household_2019}, there are around 229,000 mobility of care trips taken on public transportation every day in the Chicago region. Most of the riders on these trips are women (59\%) and low-income (66\%). By aggregating the POLARIS activity cancellation rate in several categories – errands, healthcare, pickup-dropoff, school, and shop-major – we estimate a 12\% drop (27,000 trips) in mobility of care trips without transit service. 

\subsection{Economic Impact}

We consider three illustrative categories of economic impact: 1) Consumer spending losses due to activity cancellations, 2) Travel time losses due to congestion, and 3) Annual car ownership increase. 

Daily reduction in certain activity categories were translated into annual reduction in spendings using the Bureau of Labor Statistics (BLS) data on Chicago-region consumer expenditures \citep{noauthor_consumer_nodate}. Activity categories from POLARIS were translated to BLS spending categories as follows: shop-major, shop-other, and errands were classified as apparel and services, leisure was classified as entertainment, and eat out was classified as food away from home. After multiplying annual household spending in each category by the 3.8 million households in the Chicago region, the corresponding percent reduction in activities in that category from the POLARIS results was applied to the regional spending to estimate expected losses in spending due to activity cancellation. 

The regional increase of highway modes is 1.3 million vehicle hours per day with an average occupancy of 1.48. Assuming a value of time of \$30 per hour and only using 261 annual weekdays, the annual loss is calculated at \$15.1 billion. The removal of transit results in a removal of 0.7 million passenger hours. Assuming a value of time of \$25.50 per hour – value of time riding transit is generally lower compared to value of time in a car – and only using 261 annual weekdays, the annual gain is calculated at \$4.7 billion. The net annual loss due to travel time increases is \$10.4 billion. 

According to the American Automobile Association, annual cost of owning a car in 2022 was \$10,728, which became \$12,182 in 2023 \citep{moye_annual_2023}. Using the more conservative value of \$10,728 and assuming in the transit removal scenario there will be 1.9 million additional cars, the total amount is \$20.4 billion. The grand total is \$35.4 billion as shown in \cref{economic_loss}. When we compare it to the region’s \$2.7 billion in annual transit operating funding (including federal relief) \citep{noauthor_regional_nodate}, we find that for every dollar invested in transit there is a return of \$13 in the selected economic activities and travel time savings.

\begin{table}[!ht]
  \footnotesize
  \caption{Economic losses in the region due to transit removal.}\label[tab]{economic_loss}
  \begin{center}
  \begin{tabular*}\textwidth{l@{\extracolsep{\fill}}r}
      \toprule
Category                                                     & Annual total (\$ in billion) \\
\midrule
Due to activity cancellations                                & 4.6                          \\
\qquad \textit{Entertainment} & 2.4                          \\
\qquad \textit{Restaurants}  & 1.6                          \\
\qquad \textit{Sales}         & 0.6                          \\
Travel time losses                                           & 10.4                         \\
Annual car ownership increase                                & 20.4                         \\
\midrule
Grand Total                                                  & 35.4                                     \\
      \bottomrule
    \end{tabular*}
  \end{center}
\end{table}

\section{Conclusion}
Transit is the backbone of urban transportation and is a key resource for achieving net-zero targets. A complete removal of transit does not only lead to increased car ownership and higher traffic congestion, but it also causes many activity cancellations and a shrinkage in the economy. Transit, especially in urban areas such as the Chicago Metropolitan Region considered in this study, plays a key role in enabling access, furthering transportation equity, and sustaining a dynamic and active economy.  

Since the 2020 pandemic, transit has lost 30.2\% of its riders in the Chicago region. Bringing down the transit ridership to zero by removing it entirely would result in an increase of average travel times by 14.2\% in the region and 34.7\% in the city. The removal of transit and the subsequent increase in congestion also result in major activity cancellations. In the Chicago Metropolitan Region, 11.8\% of non-work and 2.8\% of work/school activities get cancelled, resulting in an overall activity cancellation of 8.6\%. The severity of activity cancellations is exacerbated in the city with 26.9\% of non-work and 7.3\% of work/school activities getting cancelled, leading to an overall activity cancellation of 19.9\%. 

Results show that removing transit impacts different groups at varying levels. For instance, women and lower income group individuals are prone to canceling their activities more than men and higher income groups. Women account for 53.7\% of non-work activity cancellations and 53.0\% of overall activity cancellations. A large percent (50.2\%) of non-work activity cancellations is experienced by the lower income groups (the lowest 40\%). Similarly, 48.0\% of overall activity cancellations are experienced by them. Moreover, 12.0\% percent of mobility of care (unpaid labor related to caring for and keeping up a household) trips are cancelled, furthering the disproportionate effects on women and low-income residents. 

Activity cancellations, travel time losses, and car ownership increase together cost \$35.4 billion in the region. The total annual public funding for transit is \$2.7 billion, resulting in a ratio of 13 to 1, highlighting the importance of transit not only from a mobility and equity perspective but also an economic one. 

This study did not simulate the potential gradual evolution towards the transit removal through the negative feedback loop of revenue losses and service reductions. On the other hand, it directly demonstrated the impact of a such final outcome in one step. Future studies will show this gradual development in an iterative framework. Moreover, the current version or the iterative version can also be applied to other metropolitan regions in the country and in the world, facing similar challenges. The findings from this study could be paired with traditional economic impact analysis to consider additional benefits of transit. Another focus area is improving the transit service including but not limited to bus speed improvements, bus rapid transit deployment with lane removals from passenger car traffic, frequency improvement with a focus on off-peak periods, as well as major route expansion scenarios. Finally, these transit improvements can be coupled with longer-term changes in land use (transit-oriented development and urbanization of the population) as well as some households reducing their car ownership. 

\section*{Acknowledgments}
This work was supported by the U.S. Department of Energy (DOE) Vehicle Technologies Office (VTO) under the Transportation Systems and Mobility Tools Core Maintenance/Pathways to Net-Zero Regional Mobility, an initiative of the Energy Efficient Mobility Systems (EEMS) Program. The submitted manuscript has been created by [University of Washington] and the UChicago Argonne, LLC, Operator of Argonne National Laboratory (Argonne). Argonne, a U.S. Department of Energy Office of Science laboratory, is operated under Contract No. DE-AC02-06CH11357. The U.S. Government retains for itself, and others acting on its behalf, a paid-up nonexclusive, irrevocable worldwide license in said article to reproduce, prepare derivative works, distribute copies to the public, and perform publicly and display publicly, by or on behalf of the Government. 

The work has greatly benefited from the feedback and input by Peter Fahrenwald, Peter Kersten, and Amy Hofstra from the Regional Transit Authority; Molly Poppe, Cara Bader and Tom McKone from the Chicago Transit Authority; and Gabriel Barrett, Delphine Protopapas, and Anson So at the Massachusetts Institute of Technology Urban Mobility Lab. 

\section*{Author Contributions}
The authors confirm contribution to the paper as follows: agent-based network simulation development: Omer Verbas, Taner Cokyasar, Aymeric Rousseau, Joshua Auld; study conception: all authors; data collection: Scott Wainwright, Maeve Coates; designing simulation input and running simulations: Omer Verbas, Taner Cokyasar; analysis and interpretation of results: all authors; draft manuscript preparation: Omer Verbas, Taner Cokyasar, Seamus Joyce-Johnson, Scott Wainwright, Maeve Coates. All authors reviewed the results and approved the final version of the manuscript.
    



\clearpage
\bibliographystyle{elsarticle-harv} 
\bibliography{our_bib}

\vfill
\framebox{\parbox{.90\linewidth}{\scriptsize The submitted manuscript has been created by
        UChicago Argonne, LLC, Operator of Argonne National Laboratory (``Argonne'').
        Argonne, a U.S.\ Department of Energy Office of Science laboratory, is operated
        under Contract No.\ DE-AC02-06CH11357.  The U.S.\ Government retains for itself,
        and others acting on its behalf, a paid-up nonexclusive, irrevocable worldwide
        license in said article to reproduce, prepare derivative works, distribute
        copies to the public, and perform publicly and display publicly, by or on
        behalf of the Government.  The Department of Energy will provide public access
        to these results of federally sponsored research in accordance with the DOE
        Public Access Plan \url{http://energy.gov/downloads/doe-public-access-plan}.}}
\end{document}